\documentstyle[11pt,amstex,amssymb]{amsart}
\newtheorem{thm}{Theorem}[section]
\newtheorem{prop}[thm]{Proposition}
\newtheorem{lemma}[thm]{Lemma}
\newtheorem{cor}[thm]{Corollary}
\newtheorem{remark}[thm]{Remark}

\newcommand{\proof}{\noindent{\it Proof.}\enspace}

\def\bC{{\Bbb C}}
\def\bR{{\Bbb R}}
\def\bN{{\Bbb N}}
\def\bZ{{\Bbb Z}}
\def\bT{{\Bbb T}}
\def\cH{{\cal H}}
\def\cM{{\cal M}}
\def\cD{{\cal D}}
\def\cA{{\cal A}}
\def\cP{{\cal P}}
\def\cI{{\cal I}}
\def\MN{M_N({\Bbb C})}
\def\tr{{\rm tr}}
\def\rank{{\rm rank}}
\def\Det{{\rm Det}}
\def\Re{{\rm Re}\,}
\def\Im{{\rm Im}\,}
\def\eps{\varepsilon}
\def\Diag{{\rm Diag}}
\def\b1{{\bf 1}}
\def\ffi{\varphi}
\def\im{{\rm i}\,}
\begin{document}

\title[Large deviations for functions of two random projection matrices]
{{\large Large deviations for functions of two \\
\bigskip random projection matrices}}

\author[F. Hiai]{Fumio Hiai$\,^{1,2}$}
\address{Graduate School of Information Sciences,
Tohoku University, Aoba-ku, Sendai 980-8579, Japan}
\author[D. Petz]{D\'enes Petz$\,^{1,3}$}
\address{Alfr\'ed R\'enyi Institute of Mathematics, Hungarian 
Academy of Sciences, H-1053 Budapest, Re\'altanoda u. 13-15, Hungary}

\thanks{$^1\,$Supported in part by Japan-Hungary JSPS-HAS Joint Project.}

\thanks{$^2\,$Supported in part by Strategic Information and Communications
R\&D Promotion Scheme of MPHPT}

\thanks{$^3\,$Supported in part by OTKA T032662.}

\begin{abstract}
In this paper two independent and unitarily invariant projection matrices
$P(N)$ and $Q(N)$ are considered and the large deviation is proven for
the eigenvalue density of all polynomials of them as the matrix size $N$
converges to infinity. The result is formulated on  the tracial state
space $TS(\cA)$ of the universal $C^*$-algebra $\cA$ generated by two
selfadjoint projections. The random  pair $(P(N), Q(N))$ determines a
random tracial state $\tau_N \in TS(\cA)$ and $\tau_N$ satisfies the
large deviation. The rate function is in close connection with
Voiculescu's free entropy defined for pairs of projections.
\end{abstract}

\maketitle

\section*{Introduction}

Large deviation results for the empirical eigenvalue density of random
matrices started with the paper of Ben Arous and Guionnet \cite{BG} in
which generalized Wigner theorem concerning Gaussian symmetric (or 
selfadjoint) matrices was proven. The paper was followed by large
deviation results for several other kind of random matrices (as
Wishart, etc); see the monograph \cite{HP} for a detailed discussion
and the survey \cite{Gu} for more recent developments.

Up to now the typical large deviation results on random matrices have
dealt with the empirical eigenvalue density of a certain sequence of
matrices; occasionally these matrices were algebraically expressed from
two (as in \cite{PR}). In this paper two independent projection
matrices are considered and the large deviation is proven for all
polynomials (even for more general functions) of them. More precisely,
the main result is a $C^*$-algebraic formulation of large deviations for
the sequence of two random  selfadjoint projection matrices  $P(N)$ and
$Q(N)$ having independent and unitarily invariant distribution provided
moreover $\alpha:=\lim_N \rank(P(N))/N$ and $\beta:=\lim_N \rank(Q(N))/N$
exist. The main theorem is formulated on the tracial state space
$TS(\cA)$ of the universal $C^*$-algebra $\cA:=C^*(\bZ \star\bZ)$
generated by two selfadjoint projections $e$ and $f$. The random pair
$(P(N), Q(N))$  determines a random tracial state $\tau_N \in TS(\cA)$ as
follows:
$$
\tau_N (h)= \frac{1}{N}{\rm Tr}(\psi(h)),\qquad h \in \cA,
$$
where $\psi : \cA \to \MN$ is the unique $*$-homomorphism such that
$\psi(e)=P(N)$ and $\psi(f)=Q(N)$. The random $\tau_N$ induces a measure 
$\nu_N$ on $TS(\cA)$ and the sequence $\nu_N$ satisfies the large deviation 
principle in the scale $1/N^2$ with a rate function $\cI : TS(\cA) \to 
[0,\infty]$ in the ordinary sense. It is very remarkable that the rate 
function $\cI$ is in close relation with Voiculescu's free entropy 
$\chi(p,q)$ defined for a pair of projections in a $W^*$-probability space.
Namely, the GNS-construction from $(\cA,\tau)$ yields a $W^*$-probability
space $(\pi_\tau(\cA)'',\tilde\tau)$ and for the projections
$p=\pi_\tau(e)$ and $q=\pi_\tau(f)$, we have $\cI(\tau)= -\chi(p,q)$.

The result includes a bunch of traditional large deviation results for
the eigenvalue density of different polynomials of  $P(N)$ and $Q(N)$.
The corresponding rate function can be obtained from $\cI$ by the contraction
principle and computed explicitly in some examples as $P(N)Q(N)+ Q(N)P(N)$
and $aP(N)+bQ(N)$. 

The paper is organized as follows. First we establish a large deviation 
theorem for the empirical eigenvalue density of the random matrix 
$P(N)Q(N)P(N)$. This result is obtained via the joint eigenvalue density 
and the cases $\alpha+\beta$ $\le 1$ and $\ge 1$ are somewhat separated 
but treated parallel. A few facts about the Jacobi ensemble are used here. 
Since polynomials of two projections are easily controlled by the powers of 
$P(N)Q(N)P(N)$, we can move to the $C^*$-algebraic formulation mentioned
above. The tracial state space $TS(\cA)$ has a convenient representation
in terms of four numbers and a measure on $(0,1)$. The large deviation
theorem or more precisely the rate function is first identified in terms
of the representation of tracial states and the description \`a la
Voiculescu  comes afterwards. The last section is the application of the
contraction principle and contains very concrete computations.

\bigskip
\section{Joint distribution of two projections}
\setcounter{equation}{0}

Let $\MN$ be the algebra of $N\times N$ complex matrices. By an {\it $N\times N$
random projection matrix} $P$ we always mean a random orthogonal (or selfadjoint)
projection matrix, and the {\it unitary invariance} of $P$ means that the
distribution of $VPV^*$ is equal to that of $P$ for any unitary $V \in \MN$.

The aim of this section is to analyze the joint distribution of two independent
and unitarily invariant random projection matrices $P,Q$ in $\MN$, when their
ranks $\rank (P)=k$ and $\rank (Q)=l$ are fixed; we may assume that
$0\leq k\leq l\leq N$. Throughout this section, we keep these assumptions on $P$
and $Q$.
 
The joint eigenvalue distribution of $PQP$ is related to the Jacobi ensemble.
Let $(A,B)$ be an independent pair of $N\times N$ complex {\it Wishart matrices}
of $p$ degrees of freedom and of $q$ degrees of freedom, respectively, that is, 
$A=YY^*$ and $B=ZZ^*$ with complex $N\times p$ and $N\times q$ random matrices 
$Y$ and $Z$ such that $\Re Y_{ij}$, $\Im Y_{ij}$, $\Re Z_{ij}$ and $\Im Z_{ij}$ 
are independent standard Gaussians. Assume here that $p,q\geq N$. Then the 
random positive semidefinite matrix
$$
(A+B)^{-1/2}A(A+B)^{-1/2}
$$
is called an $N\times N$ {\it Jacobi ensemble} of parameter $(p-N,q-N)$. It
has the probability distribution
\begin{equation}\label{F-1.1}
{\rm Constant}\times \Det(X)^{p-N}\Det(I-X)^{q-N}
\b1_{\{0\leq X\leq I\}}(X)\,dX
\end{equation}
on the space of $N\times N$ selfadjoint matrices (see \cite[Lemma 2.1]{Col}), 
where $\b1_{\{0\leq X\leq I\}}$ denotes the characteristic function of
$\{X\in \MN:0\leq X\leq I\}$. The density formula (\ref{F-1.1}) implies
the joint distribution of the eigenvalues
$$
{\rm Constant}\times \prod_{i=1}^Nx_i^{p-N}(1-x_i)^{q-N}
\prod_{1\leq i<j\leq N}(x_i-x_j)^2\prod_{i=1}^N\b1_{[0,1]}(x_i)\,dx_i,
$$
see also \cite{Con} or \cite[Chapter 2]{Fo}.

The next lemma is from \cite[Theorem 2.2]{Col}.

\begin{lemma}\label{L-1.1}
Assume that $k+l\leq N$. Then $PQP$, when considered as a random matrix in
$M_k(\bC)=P\MN P$, has the distribution of a Jacobi ensemble of parameter
$(l-k,N-k-l)$. Hence, the joint eigenvalue distribution of the nonzero
eigenvalues of $PQP$ is given by
\begin{equation}\label{F-1.2}
{1\over Z_{N,k,l}}\prod_{i=1}^kx_i^{l-k}(1-x_i)^{N-k-l}
\prod_{1\leq i<j\leq k}(x_i-x_j)^2\prod_{i=1}^k\b1_{[0,1]}(x_i)\,dx_i
\end{equation}
with a normalization constant $Z_{N,k,l}$.
\end{lemma}

Let $(A,B)$ and $(A',B')$ be pairs of selfadjoint $N\times N$ random matrices.
We say that they have the same {\it joint distribution} if 
$$
\tr_N(h(A,B))=\tr_N(h(A',B')) \quad \hbox{almost surely}
$$
for any polynomial $h$ of two non-commuting variables, where $\tr_N$ denotes
the normalized trace on $\MN$.

Our strategy is to modify the pair $(P,Q)$ of projections in such a way
that they are easy to handle but their joint distribution does not
change. As the first step, we may assume that $(P,Q)$ are of the forms
$$
P=I_k\oplus0_{N-k},\quad Q=U(I_l\oplus0_{N-l})U^*,
$$
where $I_k\oplus0_{N-k}$ stands for the diagonal matrix whose $k$ first
diagonal entries are $1$ and the remaining are $0$, and $U$ is an $N\times N$ 
Haar-distributed random unitary matrix. In this way,  randomness belongs to 
only $Q$, while $P$ is a constant projection matrix. 

\begin{prop}\label{P-1.2} \
\begin{itemize}
\item[(a)] If $k+l\leq N$, then the joint distribution of $(P,Q)$ coincides 
with that of the pair
$$
P\quad\text{and}\quad
\bmatrix X&\sqrt{X(I_k-X)}&0&0\\
\sqrt{X(I_k-X)}&I_k-X&0&0\\0&0&I_{l-k}&0\\0&0&0&0_{N-k-l}
\endbmatrix,
$$
where $X:=\Diag(x_1,\dots,x_k)$ and $(x_1,\dots,x_k)\in [0,1]^k$ is
distributed under the distribution \eqref{F-1.2}.
\item[(b)] If $k+l>N$, then the joint distribution of $(P,Q)$ coincides  with 
that of the pair
$$
P\quad\text{and}\quad
\bmatrix I_{k+l-N}&0&0&0\\0&X&\sqrt{X(I_{N-l}-X)}&0\\
0&\sqrt{X(I_{N-l}-X)}&I_{N-l}-X&0\\0&0&0&I_{l-k}
\endbmatrix,
$$
where
$X:=\Diag(x_1,\dots,x_{N-l})$ and $(x_1,\dots,x_{N-l})$ in $[0,1]^{N-l}$ is
distributed under
\begin{equation}\label{F-1.3}
{1\over Z_{N,k,l}}
\prod_{i=1}^{N-l}x_i^{l-k}(1-x_i)^{k+l-N}
\prod_{1\leq i<j\leq N-l}(x_i-x_j)^2\prod_{i=1}^{N-l}\b1_{[0,1]}(x_i)\,dx_i.
\end{equation}
\end{itemize}
\end{prop}

\proof
(a)\enspace Assume $k+l\leq N$. By the structure theorem of two projections
(see \cite[pp.~306--308]{Ta}), after a (random) unitary conjugation, $(P,Q)$ can
be represented as
\begin{eqnarray*}
P&=&\bmatrix I&0\\0&0\endbmatrix\oplus I\oplus I\oplus0\oplus0, \\
Q&=&\bmatrix X&\sqrt{X(I-X)}\\\sqrt{X(I-X)}&I-X\endbmatrix
\oplus I\oplus0\oplus I\oplus0,
\end{eqnarray*}
where $0\leq X\leq I$ with $\ker X=\{0\}$ and $\ker(I-X)=\{0\}$ on $\cH_0$, under
a decomposition
$$
\bC^N=(\cH_0\otimes\bC^2)\oplus\cH_1\oplus\cH_2\oplus\cH_3\oplus\cH_4.
$$
(Note that $\cH_1$, $\cH_2$, $\cH_3$ and $\cH_4$ are the ranges of $P\wedge Q$,
$P\wedge Q^\perp$, $P^\perp\wedge Q$ and $(P\vee Q)^\perp$, respectively, and
some of them may be zero spaces.) Since $PQP|_{P\bC^N}$ is $X\oplus I\oplus0$ on
$\cH_0\oplus\cH_1\oplus\cH_2$, it follows from Lemma \ref{L-1.1} that $\cH_1$ 
and $\cH_2$ are zero spaces almost surely. This shows that there exists an
$N\times N$ random unitary matrix $V$ such that
\begin{eqnarray*}
VPV^*&=&P=\bmatrix I_k&0\\0&0\endbmatrix\oplus0_{l-k}\oplus0_{N-k-l}, \\
VQV^*&=&\bmatrix X&\sqrt{X(I_k-X)}\\\sqrt{X(I_k-X)}&I_k-X\endbmatrix
\oplus I_{l-k}\oplus0_{N-k-l},
\end{eqnarray*}
where $X=\Diag(x_1,\dots,x_k)$ and $(x_1,\dots,x_k)\in[0,1]^k$ is distributed
under \eqref{F-1.2} by Lemma \ref{L-1.1}. Hence we have the desired conclusion.

(b)\enspace Next, assume $k+l>N$; then since $N-l<k$ and $(N-l)+k\leq N$, one
can apply the above case (a) to $(I-Q,P)$ instead of $(P,Q)$. Thus, the joint
distribution of $(I-Q,P)$ is almost surely equal to that of the pair
$$
\bmatrix I_{N-l}&0\\0&0\endbmatrix\oplus0_{k+l-N}\oplus0_{l-k}
$$
and
$$
\bmatrix X&\sqrt{X(I_{N-l}-X)}\\\sqrt{X(I_{N-l}-X)}&I_{N-l}-X\endbmatrix
\oplus I_{n+m-N}\oplus0_{m-n}
$$
so that $(P,Q)$ has the same joint distribution almost surely as the pair
$$
\bmatrix X&\sqrt{X(I_{N-l}-X)}\\\sqrt{X(I_{N-l}-X)}&I_{N-l}-X\endbmatrix
\oplus I_{k+l-N}\oplus0_{l-k}
$$
and
$$
\bmatrix0&0\\0&I_{N-l}\endbmatrix\oplus I_{k+l-N}\oplus I_{l-k}.
$$
Here, $X=\Diag(x_1,\dots,x_{N-l})$ and $(x_1,\dots,x_{N-l})\in[0,1]^{N-l}$ is
distributed under
\begin{equation}\label{F-1.4}
{1\over Z_{N,N-l,k}}\prod_{i=1}^{N-l}x_i^{k+l-N}(1-x_i)^{l-k}
\prod_{1\leq i<j\leq N-l}(x_i-x_j)^2
\prod_{i=1}^{N-l}\b1_{[0,1]}(x_i)\,dx_i.
\end{equation}
Since $\bmatrix X&\sqrt{X(I_{N-l}-X)}\\\sqrt{X(I_{N-l}-X)}&I_{N-l}-X\endbmatrix$
and $\bmatrix0&0\\0&I_{N-l}\endbmatrix$ are respectively transformed into
$\bmatrix I_{N-l}&0\\0&0\endbmatrix$ and
$\bmatrix I_{N-l}-X&\sqrt{X(I_{N-l}-X)}\\\sqrt{X(I_{N-l}-X)}&X\endbmatrix$ by a
conjugation by the unitary matrix $\bmatrix\sqrt
X&\sqrt{I_{N-l}-X}\\-\sqrt{I_{N-l}-X}&\sqrt X\endbmatrix$, the conclusion follows
after the coordinate change $X\mapsto I_{N-l}-X$ so that \eqref{F-1.4} is
transformed into \eqref{F-1.3}.\qed

\bigskip
From Proposition \ref{P-1.2} we can readily obtain joint eigenvalue distributions
of some polynomials of $P$ and $Q$. For example, we have:

\begin{cor}\label{C-1.3} \
\begin{itemize}
\item[(i-a)] When $k+l\leq N$, the eigenvalues of $PQP$ (or $PQ$) are given as
$$
\underbrace{0,\dots,0}_{N-k\ {\rm times}},
x_1,\dots,x_k
$$
and the joint distribution of $(x_1,\dots,x_k)$ is \eqref{F-1.2}.
\item[(i-b)] When $k+l>N$, the eigenvalues of $PQP$ (or $PQ$) are given as
$$
\underbrace{0,\dots,0}_{N-k\ {\rm times}},
\underbrace{1,\dots,1}_{k+l-N\ {\rm times}},
x_1,\dots,x_{N-l},
$$
and the joint distribution of $(x_1,\dots,x_{N-l})$ is \eqref{F-1.3}.
\item[(ii-a)] When $k+l\leq N$, the eigenvalues of $PQ+QP$ are given as
$$
\underbrace{0,\dots,0}_{N-2k\ {\rm times}},
x_1\pm\sqrt{x_1},\dots,x_k\pm\sqrt{x_k}
$$
and the joint distribution of $(x_1,\dots,x_k)$ is \eqref{F-1.2}.
\item[(ii-b)] When $k+l>N$, the eigenvalues of $PQ+QP$ are given as
$$
\underbrace{0,\dots,0}_{l-k\ {\rm times}},
\underbrace{2,\dots,2}_{k+l-N\ {\rm times}},
x_1\pm\sqrt{x_1},\dots,x_{N-l}\pm\sqrt{x_{N-l}},
$$
and the joint distribution of $(x_1,\dots,x_{N-l})$ is \eqref{F-1.3}.
\item[(iii-a)] When $k+l\leq N$ and $a,b\in\bR\setminus\{0\}$, the eigenvalues of
$aP+bQ$ are given as
$$
\underbrace{0,\dots,0}_{N-2k\ {\rm times}},
\underbrace{b,\dots,b}_{l-k\ {\rm times}},
x_1,\dots,x_k,a+b-x_1,\dots,a+b-x_k,
$$
and the joint distribution of $(x_1,\dots,x_k)$ is
\begin{eqnarray}\label{F-1.5}
&&{2^k\over|ab|^{k(N-k)}Z_{N,k,l}}
\prod_{i=1}^k\bigg|x_i-{a+b\over2}\bigg|\,\big|(x_i-a)(x_i-b)\big|^{l-k}
\,\big|x_i(a+b-x_i)\big|^{N-k-l} \nonumber\\
&&\hskip3cm\times\prod_{1\leq i<j\leq k}(x_i-x_j)^2(a+b-x_i-x_j)^2
\prod_{i=1}^k\b1_{[A,B]}(x_i)\,dx_i,
\end{eqnarray}
where $Z_{N,k,l}$ is the normalization constant in \eqref{F-1.2} and $A,B$ are
the first two smallest numbers of $0,a,b,a+b$.
\item[(iii-b)] When $k+l>N$ and $a,b\in\bR\setminus\{0\}$, the eigenvalues of
$aP+bQ$ are given as
$$
\underbrace{b,\dots,b}_{l-k\ {\rm times}},
\underbrace{a+b,\dots,a+b}_{k+l-N\ {\rm times}},
x_1,\dots,x_{N-l},a+b-x_1,\dots,a+b-x_{N-l},
$$
and the joint distribution of $(x_1,\dots,x_{N-l})$ is
\begin{eqnarray*}
&&{2^{N-l}\over|ab|^{l(N-l)}Z_{N,k,l}}
\prod_{i=1}^{N-l}\bigg|x_i-{a+b\over2}\bigg|\,\big|(x_i-a)(x_i-b)\big|^{l-k}
\,\big|x_i(a+b-x_i)\big|^{k+l-N} \\
&&\qquad\qquad\qquad\times\prod_{1\leq i<j\leq N-l}(x_i-x_j)^2(a+b-x_i-x_j)^2
\prod_{i=1}^{N-l}\b1_{[A,B]}(x_i)\,dx_i,
\end{eqnarray*}
where $Z_{N,k,l}$ is the normalization constant in \eqref{F-1.3} and $A,B$ are as
in {\rm(iii-a)}.
\end{itemize}
\end{cor}

\proof
(i-a) is Lemma \ref{L-1.1} and (i-b) is immediate from Proposition
\ref{P-1.2}\,(b).

(ii-a)\enspace By Proposition \ref{P-1.2}\,(a) we may assume that
$$
PQ+QP=\bmatrix2X&\sqrt{X(I_k-X)}\\
\sqrt{X(I_k-X)}&0\endbmatrix\oplus0_{N-2k},
$$
where $X$ is as in Proposition \ref{P-1.2}\,(a). Then the result immediately
follows because the eigenvalues of the $2\times2$ matrix
$\bmatrix2x&\sqrt{x(1-x)}\\\sqrt{x(1-x)}&0\endbmatrix$ for $0\leq x\leq1$ are
$x\pm\sqrt x$. The proof of (ii-b) is similar by Proposition \ref{P-1.2}\,(b).

(iii-a)\enspace By Proposition \ref{P-1.2}\,(a) we may assume that
$$
aP+bQ=\bmatrix aI_k+bX&b\sqrt{X(I_k-X)}&0&0\\b\sqrt{X(I_k-X)}&b(I_k-X)&0&0
\\0&0&bI_{l-k}&0\\0&0&0&0_{N-k-l}
\endbmatrix.
$$
The eigenvalues of the $2\times2$ matrix
$\bmatrix a+bx&b\sqrt{x(1-x)}\\b\sqrt{x(1-x)}&b(1-x)\endbmatrix$ for $0\leq x\leq1$
are
$$
{a+b\pm\sqrt{(a-b)^2+4abx}\over2}.
$$
Set $t_i:={a+b-\sqrt{(a-b)^2+4abx_i}\over2}$ for $1\leq i\leq k$.
Then the eigenvalues of $aP+bQ$ are
$$
\underbrace{0,\dots,0}_{N-2k\ {\rm times}},
\underbrace{b,\dots,b}_{l-k\ {\rm times}},
t_1,\dots,t_k,a+b-t_1,\dots,a+b-t_k,
$$
and $(t_1,\dots,t_k)$ is supported in $[A,B]^k$. By noting that
$$
x_i={(t_i-a)(t_i-b)\over ab},\quad
1-x_i={t_i(a+b-t_i)\over ab},\quad
{dx_i\over dt_i}=2\biggl(t_i-{a+b\over2}\biggr),
$$
the form \eqref{F-1.5} of the joint distribution of $(x_1,\dots,x_k)$ can be
directly computed from \eqref{F-1.2}. The proof of (iii-b) is similar.\qed

\bigskip
\section{Large deviation for $PQP$}
\setcounter{equation}{0}

From now on, for each $N\in\bN$ let $(P(N),Q(N))$ be a pair of independent and
unitarily invariant random projection matrices in $\MN$  with non-random ranks
$k(N):=\rank(P(N))$ and $l(N):=\rank(Q(N))$. Throughout what follows, we assume
that $k(N)/N\to\alpha$ and $l(N)/N\to\beta$ as $N\to\infty$ for some
$\alpha,\beta\in[0.1]$. Our goal is to  obtain a large deviation theorem for the
empirical eigenvalue density of  $P(N)Q(N)P(N)$. Concerning large deviation
theory, our general reference is \cite{DZ}, but \cite{HP} contains many matrix
examples.

We have already observed that the two cases $\alpha + \beta \le 1$ and  $\alpha + 
\beta \ge 1$ are slightly different. To treat them parallel, we set
$$
n_0(N):=N-\min\{k(N),l(N)\},\quad
n_1(N):=\max\{k(N)+l(N)-N,0\},
$$
$$
n(N):=N-n_0(N)-n_1(N)\ \bigl(=\min\{k(N),l(N),N-k(N),N-l(N)\}\bigr).
$$
Then one can combine (i-a) and (i-b) of Corollary \ref{C-1.3} to see that the
eigenvalues of the $N\times N$ selfadjoint random matrix $P(N)Q(N)P(N)$ are
$$
\underbrace{0,\dots,0}_{n_0(N)\ {\rm times}},
\underbrace{1,\dots,1}_{n_1(N)\ {\rm times}}, \, x_1,\dots,x_{n(N)}
$$
and the joint distribution of $(x_1,\dots,x_{n(N)})$ is
\begin{equation}\label{F-2.1}
{1\over Z(N)}\prod_{i=1}^{n(N)}x_i^{|k(N)-l(N)|}(1-x_i)^{|k(N)+l(N)-N|}
\prod_{1\leq i<j\leq n(N)}(x_i-x_j)^2\prod_{i=1}^{n(N)}\b1_{[0,1]}(x_i)\,dx_i
\end{equation}
with a normalization constant $Z(N)$.

When ${\cal X}$ is a Polish space, let $\cM({\cal X})$ denote the set of all
probability measures on ${\cal X}$, which becomes a Polish space with respect
to weak topology. For $\mu\in\cM(\bR)$ let $\Sigma(\mu)$ be Voiculescu's
{\it free entropy} (or the minus of the {\it logarithmic energy}) of $\mu$
defined by
$$
\Sigma(\mu):=\iint\log|x-y|\,d\mu(x)\,d\mu(y)
$$
(see \cite{V1} and \cite[\S5.3]{HP}). In particular, when $\mu$ is compactly
supported, $\Sigma(\mu)\in[-\infty,+\infty)$ is well defined.

We first prove large deviation for the sequence of distributions (\ref{F-2.1})
with slight modifications of notation.

\begin{prop}\label{P-2.1}
For each $N\in\bN$ consider the distribution
\begin{equation}\label{F-2.2}
{1\over Z(N)}\prod_{i=1}^{n(N)}x_i^{\kappa(N)}(1-x_i)^{\lambda(N)}
\prod_{1\leq i<j\leq n(N)}(x_i-x_j)^2\prod_{i=1}^{n(N)}\b1_{[0,1]}(x_i)\,dx_i
\end{equation}
on $[0,1]^{n(N)}$ with $n(N)\in\bN$, $\kappa(N),\lambda(N)\in[0,\infty)$ and a
normalization constant $Z(N)$. Assume that $n(N)/N\to\rho$,
$\kappa(N)/N\to\kappa$ and $\lambda(N)/N\to\lambda$ as $N\to\infty$ for some
$\rho\in(0,\infty)$ and $\kappa,\lambda\in[0,\infty)$. Then:
\begin{itemize}
\item[(1)] The limit $\lim_{N\to\infty}{1\over N^2}\log Z(N)$ exists and it
equals $\rho^2B(\kappa/\rho,\lambda/\rho)$, where
\begin{eqnarray*}
B(s,t)&:=&{(1+s)^2\over2}\log(1+s)-{s^2\over2}\log s
+{(1+t)^2\over2}\log(1+t)-{t^2\over2}\log t \\
&&\quad-{(2+s+t)^2\over2}\log(2+s+t)+{(1+s+t)^2\over2}\log(1+s+t)
\end{eqnarray*}
for $s,t\geq0$.
\item[(2)] When $(x_1,\dots,x_{n(N)})$ is distributed under \eqref{F-2.2}, the
empirical measure
\begin{equation}\label{F-2.3}
{\delta_{x_1}+\dots+\delta_{x_{n(N)}}\over n(N)}
\end{equation}
satisfies the large deviation principle in the scale $1/N^2$ with the rate
function
\begin{equation}\label{F-2.4}
I(\mu):=-\rho^2\Sigma(\mu)
-\rho\int_0^1\bigl(\kappa\log x+\lambda\log(1-x)\bigr)\,d\mu(x)
+\rho^2B\biggl({\kappa\over\rho},{\lambda\over\rho}\biggr)
\end{equation}
for $\mu\in\cM([0,1])$. Moreover, there exists a unique minimizer
$\mu_0\in\cM([0,1])$ of $I(\mu)$ with $I(\mu_0)=0$.
\end{itemize}
\end{prop}

\proof
(1)\enspace The Selberg integral formula (see \cite[\S17.1]{Me}) gives
\begin{eqnarray*}
Z(N)&=&\int_{[0,1]^{n(N)}}
\prod_{i=1}^{n(N)}x_i^{\kappa(N)}(1-x_i)^{\lambda(N)}
\prod_{1\leq i<j\leq n(N)}(x_i-x_j)^2\prod_{i=1}^{n(N)}dx_i \\
&=&\prod_{j=1}^{n(N)}{\Gamma(j+1)\Gamma(j+\kappa(N))\Gamma(j+\lambda(N))
\over\Gamma(2)\Gamma(j+n(N)+\kappa(N)+\lambda(N))}.
\end{eqnarray*}
By using the Stirling formula, under neglecting the small order $o(N)$, we
compute
\begin{eqnarray*}
&&{1\over N^2}\log Z(N) \\
&&\quad={1\over N^2}\Biggl\{\sum_{j=1}^{n(N)}j\log j
+\sum_{j=1}^{n(N)}(j+\kappa(n))\log(j+\kappa(n))
+\sum_{j=1}^{n(N)}(j+\rho(n))\log(j+\rho(n)) \\
&&\hskip2cm-\sum_{j=1}^{n(N)}(j+n+\kappa(n)+\rho(n))
\log(j+n+\kappa(n)+\rho(n))\Biggr\} \\
&&\quad={n(N)\over N^2}\Biggl\{\sum_{j=1}^{n(N)}{j\over n(N)}\log{j\over n(N)}
+\sum_{j=1}^{n(N)}\biggl({j\over n(N)}+{\kappa\over\rho}\biggr)
\log\biggl({j\over n(N)}+{\kappa\over\rho}\biggr) \\
&&\hskip2cm+\sum_{j=1}^{n(N)}
\biggl({j\over n(N)}+{\lambda\over\rho}\biggr)
\log\biggl({j\over n(N)}+{\lambda\over\rho}\biggr) \\
&&\hskip2cm-\sum_{j=1}^{n(N)}
\biggl({j\over n(N)}+1+{\kappa\over\rho}+{\lambda\over\rho}\biggr)
\log\biggl({j\over n(N)}+1+{\kappa\over\rho}+{\lambda\over\rho}\biggr)
\Biggr\}.
\end{eqnarray*}
Therefore,
\begin{eqnarray*}
&&\lim_{N\to\infty}{1\over N^2}\log Z(N) \\
&&\quad=\rho^2\Biggl\{\int_0^1x\log x\,dx
+\int_0^1\biggl(x+{\kappa\over\rho}\biggr)
\log\biggl(x+{\kappa\over\rho}\biggr)\,dx
+\int_0^1\biggl(x+{\lambda\over\rho}\biggr)
\log\biggl(x+{\lambda\over\rho}\biggr)\,dx \\
&&\hskip2cm-\int_0^1\biggl(x+1+{\kappa\over\rho}+{\lambda\over\rho}\biggr)
\log\biggl(x+1+{\kappa\over\rho}+{\lambda\over\rho}\biggr)\,dx
\Biggr\} \\
&&\quad=\rho^2B\biggl({\kappa\over\rho},{\lambda\over\rho}\biggr).
\end{eqnarray*}

(2)\enspace Denote the distribution \eqref{F-2.2} by $\nu_{n(N)}$ and define the
probability measure $P_N$ on $\cM([0,1])$ by
$$
P_N(\Lambda)
:=\nu_{n(N)}\bigl(\bigl\{x\in[0,1]^{n(N)}:\mu_x\in\Lambda\bigr\}\bigr)
$$
for Borel subsets $\Lambda$ of $\cM([0,1])$, where $\mu_x$ denotes the empirical
measure \eqref{F-2.3} for $x=(x_1,\dots,x_{n(N)})$. Define the kernel functions
on $[0,1]^2$ as follows:
$$
F(x,y):=-\log|x-y|-{\kappa\over2\rho}(\log x+\log y)
-{\lambda\over2\rho}(\log(1-x)+\log(1-y)),
$$
$$
F_R(x,y):=\min\{F(x,y),R\}\quad\mbox{for $R>0$}.
$$
Furthermore, for each $N\in\bN$ we define
\begin{eqnarray*}
\widetilde F_N(x,y)&:=&-\log|x-y|
-\delta_{\kappa>0}{\kappa(N)\over2n(N)}(\log x+\log y) \\
&&\qquad-\delta_{\lambda>0}{\lambda(N)\over2n(N)}(\log(1-x)+\log(1-y)),
\end{eqnarray*}
$$
\widetilde F_{N,R}(x,y):=\min\{\widetilde F_N(x,y),R\}
\quad\mbox{for $R>0$},
$$
where $\delta_{\kappa>0}=1$ if $\kappa>0$, $\delta_{\kappa>0}=0$ if $\kappa=0$,
and $\delta_{\lambda>0}$ is similar. Then we observe the following:
\begin{itemize}
\item[(i)] $\widetilde F_{N,R}(x,y)\le-\log|x-y|-{\kappa(n)\over2n(N)}(\log
x+\log y)-{\lambda(n)\over2n(N)}(\log(1-x)+\log(1-y))$ for all $x,y\in[0,1]$.
\item[(ii)] For any $R>0$, $\widetilde F_{N,R}(x,y)$ converges to $F_R(x,y)$
uniformly for $x,y\in[0,1]$ as $N\to\infty$.
\end{itemize}
In fact, (i) is obvious by the definition of $\widetilde F_{N,R}(x,y)$. For (ii)
assume that $\kappa,\lambda>0$ (the proof is similar for other cases). For
$\delta>0$ set
$$
T_\delta:=\{(x,y)\in[0,1]^2:
\delta\le x\le1-\delta,\,\delta\le y\le1-\delta,\,|x-y|\ge\delta\}.
$$
For any $R>0$ there exist $\delta>0$ and $N_0\in\bN$ such that $F(x,y)\ge R$
and $\widetilde F_N(x,y)\ge R$ for all $(x,y)\in[0,1]^2\setminus T_\delta$ and
$N\ge N_0$. Obviously, $\widetilde F_N(x,y)$ converges to $F(x,y)$ uniformly on
$T_\delta$ as $N\to\infty$, and the assertion follows.

According to general theory of large deviations (\cite{DZ}), the stated large
deviation is shown when we prove the following two inequalities for every
$\mu\in\cM([0,1])$:
\begin{equation}\label{F-2.5}
\inf_G\biggl[\limsup_{N\to\infty}{1\over N^2}\log P_N(G)\biggr]
\le-\rho^2\iint F(x,y)\,d\mu(x)\,d\mu(y)-C,
\end{equation}
\begin{equation}\label{F-2.6}
\inf_G\biggl[\liminf_{N\to\infty}{1\over N^2}\log P_N(G)\biggr]
\ge-\rho^2\iint F(x,y)\,d\mu(x)\,d\mu(y)-C,
\end{equation}
where $C:=\rho^2B(\kappa/\rho,\lambda/\rho)$ and $G$ runs over neighborhoods of
$\mu$.

\medskip
\noindent{\it Proof of \eqref{F-2.5}.}\enspace
For every neighborhood $G$ of $\mu\in\cM([0,1])$, setting $\widetilde
G:=\{x\in[0,1]^{n(N)}:\mu_x\in G\}$, by the above (i) we have
\begin{eqnarray*}
&&P_N(G)=\nu_{n(N)}(\widetilde G) \\
&&\quad={1\over Z(N)}\int_{\widetilde G}
\prod_{i=1}^{n(N)}x_i^{\kappa(N)}(1-x_i)^{\lambda(N)}
\prod_{1\le i<j\le n(N)}(x_i-x_j)^2\prod_{i=1}^{n(N)}dx_i \\
&&\quad\le{1\over Z(N)}\int_{\widetilde G}\prod_{i=1}^{n(N)}
x_i^{\kappa(N)/n(N)}(1-x_i)^{\lambda(N)/n(N)} \\
&&\hskip3cm\times
\exp\Biggl(-2\sum_{1\le i<j\le n(N)}\widetilde F_{N,R}(x_i,x_j)
\Biggr)\prod_{i=1}^{n(N)}dx_i \\
&&\quad\le{1\over Z(N)}\Biggl(\int_0^1
x^{\kappa(N)/n(N)}(1-x)^{\lambda(N)/n(N)}\,dx\Biggr)^{n(N)} \\
&&\hskip1cm\times\exp\Biggl(-n(N)^2\inf_{\mu'\in G}\iint
\widetilde F_{N,R}(x,y)\,d\mu'(x)\,d\mu'(y)+n(N)R\Biggr).
\end{eqnarray*}
Since the above fact (ii) implies that
$$
\lim_{N\to\infty}\Biggl(\inf_{\mu'\in G}
\iint\widetilde F_{N,R}(x,y)\,d\mu'(x)\,d\mu'(y)\Biggr)
=\inf_{\mu'\in G}\iint F_R(x,y)\,d\mu'(x)\,d\mu'(y),
$$
we get
$$
\lim_{N\to\infty}{1\over N^2}\log P_N(G)
\le-\rho^2\inf_{\mu'\in G}\iint F_R(x,y)\,d\mu'(x)\,d\mu'(y)-C
$$
thanks to (1). Furthermore, appealing to the continuity of
$\mu'\mapsto\iint F_R(x,y)\,d\mu'(x)\,d\mu'(y)$, we obtain
$$
\inf_G\biggl[\limsup_{N\to\infty}{1\over N^2}\log P_N(G)\biggr]
\le-\rho^2\iint F_R(x,y)\,d\mu(x)\,d\mu(y)-C
$$
so that \eqref{F-2.5} follows by letting $R\to+\infty$.

\medskip
\noindent{\it Proof of \eqref{F-2.6}.}\enspace
If $\mu$ has an atom at $0$ or $1$, then
$\iint F(x,y)\,d\mu(x)\,d\mu(y)=+\infty$ so that we have nothing to do.
Otherwise, letting $d\mu_\delta(x):=
\mu([\delta,1-\delta])^{-1}\b1_{[\delta,1-\delta]}(x)\,d\mu(x)$, we get
$$
\iint F(x,y)\,d\mu(x)\,d\mu(y)
=\lim_{\delta\searrow0}\iint F(x,y)\,d\mu_\delta(x)\,d\mu_\delta(y).
$$
Also it is immediate to see that
$$
\mu\in\cM([0,1])\mapsto\inf\biggl\{\liminf_{N\to\infty}
{1\over N^2}\log P_N(G):\mbox{$G$ is a neighborhood of $\mu$}\biggr\}
$$
is upper semicontinuous. Hence we may assume that $\mu$ is supported in $[a,b]$
with $0<a<b<1$. For $\eps>0$ let $\phi_\eps\ge0$ be a $C^\infty$-function
supported in $[-\eps,\eps]$ such that $\int\phi_\eps(x)\,dx=1$. Then we get
$\Sigma(\phi_\eps*\mu)\ge\Sigma(\mu)$ (see \cite[p.~216]{HP}) as well as
$$
\lim_{\eps\searrow0}\int\log x\,d(\phi_\eps*\mu)(x)=\int\log x\,d\mu(x),
$$
$$
\lim_{\eps\searrow0}\int\log(1-x)\,d(\phi_\eps*\mu)(x)=\int\log(1-x)\,d\mu(x)
$$
so that $\mu$ may be assumed to have a continuous density. Furthermore, by the
concavity of $\Sigma(\mu)$, it suffices to prove \eqref{F-2.6} for
$(1-\eps)\mu+\eps m$ for each $0<\eps<1$, where $m$ is the uniform
measure on an interval including the support ${\rm supp}\,\mu$. After all, we
can assume that $\mu$ has a continuous density $f>0$ on ${\rm supp}\,\mu=[a,b]$
with $0<a<b<1$ and $\delta\le f(x)\le\delta^{-1}$ on $[a,b]$ for some
$\delta>0$.

For each $N\in\bN$ let
$$
a<a_1^{(N)}<b_1^{(N)}<a_2^{(N)}<\dots<a_{n(N)}^{(N)}<b_{n(N)}^{(N)}
$$
be such that
$$
\int_a^{a_i^{(N)}}f(x)\,dx={i-{1\over2}\over n(N)},\quad
\int_a^{b_i^{(N)}}f(x)\,dx={i\over n(N)},\qquad
1\le i\le n(N);
$$
then
$$
b_i^{(N)}-a_i^{(N)}\ge{\delta\over2n(N)},\qquad 1\le i\le n(N).
$$
Define
$$
\Delta_{n(N)}:=\bigl\{x=(x_1,\dots,x_{n(N)})\in[0,1]^{n(N)}:
a_i^{(N)}\le x_i\le b_i^{(N)},\,1\le i\le n(N)\bigr\}.
$$
For any neighborhood $G$ of $\mu$, whenever $N$ is large enough, we have
$$
\Delta_{n(N)}\subset\widetilde G:=\bigl\{x\in[0,1]^{n(N)}:\mu_x\in G\bigr\}
$$
so that
\begin{eqnarray*}
&&P_N(G)=\nu_{n(N)}(\widetilde G) \\
&&\quad\ge{1\over Z(N)}\int_{\Delta_{n(N)}}
\prod_{i=1}^{n(N)}x_i^{\kappa(N)}(1-x_i)^{\lambda(N)}
\prod_{1\le i<j\le n(N)}(x_i-x_j)^2\prod_{i=1}^{n(N)}dx_i \\
&&\quad\ge{1\over Z(N)}\biggl({\delta\over2n(N)}\biggr)^{n(N)}\prod_{i=1}^{n(N)}
\bigl(a_i^{(N)}\bigr)^{\kappa(N)}\bigl(1-b_i^{(N)}\bigr)^{\lambda(N)}
\prod_{1\le i<j\le n(N)}\bigl(a_j^{(N)}-b_i^{(N)}\bigr)^2.
\end{eqnarray*}
With $g:[0,1]\to[a,b]$ being the inverse function of
$t\in[a,b]\mapsto\int_a^tf(x)\,dx$, since
$a_i^{(N)}=g\bigl(\bigl(i-{1\over2}\bigr)/n(N)\bigr)$ and
$b_i^{(N)}=g\bigl(i/n(N)\bigr)$, we have
$$
\lim_{N\to\infty}{\kappa(N)\over N^2}
\sum_{i=1}^{n(N)}\log a_i^{(N)}
=\rho\kappa\int_0^1\log g(t)\,dt=\rho\kappa\int\log x\,d\mu(x),
$$
$$
\lim_{N\to\infty}{\kappa(N)\over N^2}
\sum_{i=1}^{n(N)}\log\bigl(1-b_i^{(N)}\bigr)
=\rho\kappa\int_0^1\log(1-g(t))\,dt=\rho\kappa\int\log(1-x)\,d\mu(x),
$$
\begin{eqnarray*}
&&\lim_{N\to\infty}{2\over N^2}
\sum_{1\le i<j\le n(N)}\log\bigl(a_j^{(N)}-b_i^{(N)}\bigr) \\
&&\qquad=2\rho^2\iint_{0\le s<t\le t\le1}
\log(g(t)-g(s))\,ds\,dt=\rho^2\Sigma(\mu).
\end{eqnarray*}
These estimates altogether imply \eqref{F-2.6}.

The proof of the large deviation is now completed, and the existence of a unique
minimizer of the rate function is known as a general result on weighted
logarithmic energy functionals (see \cite[I.1.3]{ST}).\qed

\bigskip
Now, the large deviation theorem for the random matrix $P(N)Q(N)P(N)$ can be
easily shown from Proposition \ref{P-2.1}. Set
\begin{equation}\label{F-2.7}
\rho:=\min\{\alpha,\beta,1-\alpha,1-\beta\},
\end{equation}
\begin{equation}\label{F-2.8}
C
:=\rho^2B\biggl({|\alpha-\beta|\over\rho},{|\alpha+\beta-1|\over\rho}\biggr)
\end{equation}
(meant zero if $\rho=0$), and denote by $\cM((0,1))$ the set of all probability
measures on $[0,1]$ with no atoms at $0$ and $1$.

\begin{thm}\label{T-2.2}
The empirical eigenvalue density of $P(N)Q(N)P(N)$ satisfies the large
deviation principle in the scale $1/N^2$ with the rate function
$\tilde I(\tilde\mu)$ for $\tilde\mu\in\cM([0,1])$ given as follows: If
$$
\tilde\mu=(1-\min\{\alpha,\beta\})\delta_0+\max\{\alpha+\beta-1,0\}\delta_1
+\rho\mu
$$ 
with $\mu\in\cM((0,1))$, then
\begin{eqnarray}\label{F-2.9}
\tilde I(\tilde\mu)&:=&-\rho^2\Sigma(\mu)
-\rho|\alpha-\beta|\int_0^1\log x\,d\mu(x) \nonumber\\
&&\qquad-\rho|\alpha+\beta-1|\int_0^1\log(1-x)\,d\mu(x)+C;
\end{eqnarray}
otherwise $\tilde I(\tilde\mu)=+\infty$. Moreover, a unique minimizer of
$\tilde I(\tilde\mu)$ is given by
\begin{equation}\label{F-2.10}
\tilde\mu_0:=(1-\min\{\alpha,\beta\})\delta_0+\max\{\alpha+\beta-1,0\}\delta_1
+{\sqrt{(x-\xi)(\eta-x)}\over2\pi x(1-x)}\b1_{(\xi,\eta)}(x)\,dx
\end{equation}
where
\begin{equation}\label{F-2.11}
\xi,\eta:=\alpha+\beta-2\alpha\beta
\pm\sqrt{4\alpha\beta(1-\alpha)(1-\beta)}.
\end{equation}
In particular, when $\rho=0$, $\tilde I(\tilde\mu)$ is identically $+\infty$
except at only
$\tilde\mu_0=(1-\min\{\alpha,\beta\})\delta_0+\max\{\alpha+\beta-1,0\}\delta_1$.
\end{thm}

\proof
From the fact mentioned at the beginning of the section, the empirical eigenvalue
density of $P(N)Q(N)P(N)$ is given by
$$
\widetilde R_N:={n_0(N)\over N}\delta_0+{n_1(N)\over N}\delta_1
+{n(N)\over N}R_N,
$$
where $R_N:={1\over n(N)}(\delta_{x_1}+\dots+\delta_{x_{n(N)}})$ and the joint
distribution of $(x_1,\dots,x_{n(N)})$ is \eqref{F-2.1}. First, assume that
$\rho>0$. Proposition \ref{P-2.1} says that $(R_N)$ satisfies the large
deviation in the scale $1/N^2$ with the rate function $I(\mu)$ for
$\mu\in\cM([0,1])$ given in \eqref{F-2.4} with $\kappa:=|\alpha-\beta|$ and
$\lambda:=|\alpha+\beta-1|$. We now proceed as in the proof of \cite[5.5.11]{HP}.
Let $P_N$ and $\widetilde P_N$ be the distributions on $\cM([0,1])$ of $R_N$ and
$\widetilde R_N$, respectively; then
$$
\widetilde P_N(\Lambda)=P_N\biggl(\biggl\{\mu\in\cM([0,1]):
{n_0(N)\over n}\delta_0+{n_1(N)\over N}\delta_1+{n(N)\over N}\mu\in\Lambda
\biggr\}\biggr)
$$
for $\Lambda\subset\cM([0,1])$. Let $\cD$ denote the set
$\{\rho_0\delta_0+\rho_1\delta_1+\rho\mu:\mu\in\cM([0,1])\}$, where
$\rho_0:=1-\min\{\alpha,\beta\}$ and $\rho_1:=\max\{\alpha+\beta-1,0\}$. If
$\tilde\mu\notin\cD$, then $\tilde\mu(\{0\})<\rho_0$ or
$\tilde\mu(\{1\})<\rho_1$ so that letting $\tilde\mu(\{0\})<\eps<\rho_0$ (or
$\tilde\mu(\{1\})<\eps<\rho_1$) we have a neighborhood $\widetilde
G:=\{\mu'\in\cM([0,1]):\mu'(\{0\})<\eps\ ({\rm or}\ \mu'(\{1\})<\eps)\}$ of
$\mu$. Since $\widetilde P_N(\widetilde G)=0$ for large $N$, we get
$\lim_{N\to\infty}{1\over N^2}\log\widetilde P_N(\widetilde G)=-\infty$. Next,
assume that $\tilde\mu\in\cD$ and
$\tilde\mu=\rho_0\delta_0+\rho_1\delta_1+\rho\mu$. For any neighborhood
$\widetilde G$ of $\tilde\mu$ there exists a neighborhood $G$ of $\mu$ such that
${n_0(N)\over N}\delta_0+{n_1(N)\over N}\delta_1+{n(N)\over N}G\subset
\widetilde G$ for large $N$ and hence
$$
\liminf_{N\to\infty}{1\over N^2}\log\widetilde P_N(\widetilde G)
\ge\liminf_{N\to\infty}{1\over N^2}\log P_N(G)\ge-I(\mu).
$$
On the other hand, for any neighborhood $G$ of $\mu$ there exists a neighborhood
$\widetilde G$ of $\tilde\mu$ such that
$$
\biggl({N\over n(N)}\widetilde G-{n_0(N)\over n(N)}\delta_0
-{n_1(N)\over n(N)}\delta_1\biggr)\cap\cM([0,1])\subset G,
$$
that is,
$$\biggl\{\mu\in\cM([0,1]):{n_0(N)\over n}\delta_0+{n_1(N)\over n}\delta_1
+{n(N)\over n}\mu\in\widetilde G\biggr\}\subset G
$$
for large $N$. Therefore,
$$
\inf_{\widetilde G}\biggl[\limsup_{N\to\infty}
{1\over n^2}\log\widetilde P_N(\widetilde G)\biggr]
\le\inf_G\biggl[\limsup_{N\to\infty}{1\over n^2}\log P_N(G)\biggr]
\le-I(\mu).
$$
Noting that $\Sigma(\mu)=-\infty$ if $\mu\in\cM([0,1])$ has an atom at $0$ or
$1$, we obtain the desired large deviation for $(\widetilde R_N)$ when $\rho>0$.
The proof in the case $\rho=0$ is similar to the above argument for
$\tilde\mu\notin\cD$.

Finally, the existence of a unique minimizer of $\tilde I(\tilde\mu)$ is already
known by Proposition \ref{P-2.1}. To obtain the explicit form of the minimizer,
we may apply a standard method in free probability theory. In fact, by the
{\it asymptotic freeness} due to Voiculescu \cite[Theorem 3.11]{Vo} (see also
\cite[4.3.5]{HP}), the joint distribution of $(P(N),Q(N))$ converges to that of
$(p,q)$ where $p$ and $q$ are free projections in a tracial $W^*$-probability
space $(\cM,\tau)$ with $\tau(p)=\alpha$ and $\tau(q)=\beta$. The computation by
use of $S$-transform in \cite{VDN} says that the measure \eqref{F-2.10} is the
distribution measure of $pqp$; hence it is the minimizer of
$\tilde I(\tilde\mu)$.\qed

\bigskip
Note that the rate function $\tilde I(\tilde\mu)$ is indeed lower semicontinuous
and convex on $\cM([0,1])$, which is of course a good rate function because of
the compactness of $\cM([0,1])$.

\bigskip
\section{$C^*$-algebra formulation}
\setcounter{equation}{0}

The two-dimensional commutative $C^*$-algebra $\bC\oplus\bC=C^*(\bZ_2)$ is the
universal $C^*$-algebra generated by a single orthogonal projection; hence the
universal $C^*$-algebra generated two orthogonal projections is
$$
(\bC\oplus\bC)\star(\bC\oplus\bC)=C^*(\bZ\star\bZ)
$$
with projection generators $(1,0)$'s in two components. As pointed out in
\cite[p.~14]{ABH}, one can see from the structure theorem for two projections
(\cite[pp.~306--308]{Ta}) that $C^*(\bZ\star\bZ)$ is isomorphic to an
algebra of $M_2(\bC)$-valued continuous functions on $[0,1]$; namely
$$
\cA:=\bigl\{a\in C([0,1];M_2(\bC)):
\mbox{$a(0)$ and $a(1)$ are diagonal}\bigr\},
$$
where the corresponding two projection generators are represented as
$$
e(t):=\bmatrix1&0\\0&0\endbmatrix,\quad
f(t):=\bmatrix t&\sqrt{t(1-t)}\\\sqrt{t(1-t)}&1-t\endbmatrix
\quad\mbox{for $0\le t\le1$}.
$$
We thus consider the above $C^*$-algebra $\cA$ with generators $e,f$ as the
universal $C^*$-algebra generated by two projections. We denote by $TS(\cA)$ the
set of all tracial states on $\cA$, which becomes a Polish space with respect to
w*-topology. The following lemma is a concrete description of $TS(\cA)$, the 
details are left to the reader.

\begin{lemma}
For each $\tau\in TS(\cA)$ there exist
$\alpha_{11},\alpha_{10},\alpha_{01},\alpha_{00}\ge0$ with
$\sum_{i,j=0}^1\alpha_{ij}\le1$ and $\mu\in\cM((0,1))$ such that
\begin{eqnarray*}
\tau(a)&=&\alpha_{10}a_1(0)+\alpha_{01}a_2(0)
+\alpha_{11}a_1(1)+\alpha_{00}a_2(1) \\
&&\qquad+\Biggl(1-\sum_{i,j=0}^1\alpha_{ij}\Biggr)\int_0^1\tr (a(t))\,d\mu(t)
\end{eqnarray*}
for all $a\in\cA$ with $a(0)=\Diag(a_1(0),a_2(0))$ and
$a(1)=\Diag(a_1(1),a_2(1))$.
\end{lemma}

In this way, the set $TS(\cA)$ is parameterized by the set of all
$(\{\alpha_{ij}\}_{i,j=0}^1,\mu)$ of $\alpha_{ij}\ge0$,
$\sum_{i,j=0}^1\alpha_{ij}\le1$ and $\mu\in\cM((0,1))$, and we write
$\tau=(\{\alpha_{ij}\}_{i,j=0}^1,\mu)$ under this parameterization. But, note that
$\mu$ is irrelevant if $\sum_{i,j=0}^1\alpha_{ij}=1$. For
$\tau=(\{\alpha_{ij}\}_{i,j=0}^1,\mu)$ we have
$$
\tau(e)={1\over2}(1+\alpha_{11}+\alpha_{10}-\alpha_{01}-\alpha_{00}),
$$
$$
\tau(f)={1\over2}(1+\alpha_{11}-\alpha_{10}+\alpha_{01}-\alpha_{00}).
$$
Furthermore, let $\pi_\tau$ be the GNS representation of $\cA$ associated with
$\tau$ and $\tilde\tau$ be the normal extension of $\tau$ to $\pi_\tau(\cA)''$.
Then, for $p:=\pi_\tau(e)$ and $q:=\pi_\tau(f)$ in $\pi_\tau(\cA)''$ we have
\begin{equation}\label{F-3.1}
\tilde\tau(p\wedge q)=\alpha_{11},\quad
\tilde\tau(p\wedge q^\perp)=\alpha_{10},\quad
\tilde\tau(p^\perp\wedge q)=\alpha_{01},\quad
\tilde\tau(p^\perp\wedge q^\perp)=\alpha_{00}.
\end{equation}

For any two projections $p,q$ in a tracial $W^*$-probability space $(\cM,\tau)$,
the universality property of $\cA$ shows that there exists a (unique)
$*$-homomorphism $\psi_{p,q}:\cA\to\cM$ such that $\psi_{p,q}(e)=p$ and
$\psi_{p,q}(f)=q$. We simply write $h(p,q)$ for $\psi_{p,q}(h)$ for each
$h\in\cA$, which may be regarded as a sort of ``noncommutative functional
calculus." Then a tracial state $\tau_{p,q}\in TS(\cA)$ is defined by
$\tau_{p,q}(h):=\tau(h(p,q))$ for $h\in\cA$. In particular, for $N\times N$
projection matrices $P,Q$, we have $\tau_{P,Q}\in TS(\cA)$ given by
$\tau_{P,Q}(h)=\tr_N(h(P,Q))$ for $h\in\cA$. When $P,Q$ are random projection
matrices, $\tau_{P,Q}$ is a random tracial state on $\cA$ regarded as the
``noncommutative empirical measure" of the pair $(P,Q)$. Its {\it distribution
measure} on $TS(\cA)$ is defined by
$$
\nu(\Lambda):={\rm Prob}(\{\tau_{P,Q}\in\Lambda\})
$$
for Borel subsets $\Lambda\subset TS(\cA)$, where ${\rm Prob}$ denotes
probability measure of the underlying probability space where $P,Q$ are defined.

We are now in a position to state our main large deviation result formulated on
the tracial state space $TS(\cA)$.

\begin{thm}\label{T-3.2}
For each $N\in\bN$ let $(P(N),Q(N))$ be a pair of independent and
unitarily invariant random projection matrices in $\MN$ such that
$\rank (P(N))/N\to\alpha$ and $\rank (Q(N))/N\to\beta$ as $N\to\infty$.
Let $\nu_N$ be the distribution measure of the random tracial state
$\tau_N:=\tau_{P(N),Q(N)}$ on $TS(\cA)$. Then $(\nu_N)$ satisfies 
the large deviation principle in the scale $1/N^2$ with rate 
function 
\begin{eqnarray*}
\cI(\tau)&:=&-\rho^2\Sigma(\mu)
-\rho|\alpha-\beta|\int_0^1\log x\,d\mu(x) \\
&&\qquad-\rho|\alpha+\beta-1|\int_0^1\log(1-x)\,d\mu(x)+C
\end{eqnarray*}
evaluated at $\tau=(\{\alpha_{ij}\}_{i,j=0}^1,\mu)\in TS(\cA)$ if
\begin{equation}\label{F-3.2}
\cases\alpha_{11}=\max\{\alpha+\beta-1,0\}, \\
\alpha_{00}=\max\{1-\alpha-\beta,0\}, \\
\alpha_{10}=\max\{\alpha-\beta,0\}, \\
\alpha_{01}=\max\{\beta-\alpha,0\},
\endcases
\end{equation}
otherwise $\cI (\tau)=+\infty$. (See \eqref{F-2.7} and \eqref{F-2.8} for
constants $\rho$ and $C$.)

Moreover, the unique minimizer of $\cI$ is the tracial state $\tau_{p,q}$ 
corresponding to a pair $(p,q)$ of free projections with trace values 
$\alpha$ and $\beta$.
\end{thm}

\proof
First we notice that all mixed moments of $e,f$ with respect to $\tau$ are
listed as $\tau(e)$, $\tau(f)$ and
\begin{equation}\label{F-3.3}
\tau((ef)^k)=\tau((fe)^k)=\tau((efe)^k)=\tau((fef)^k),\quad k\ge1.
\end{equation}
Since the moments $\tau((efe)^k)$, $k\ge1$, determine the distribution of $efe$
with respect to $\tau$, one can define an affine homeomorphism $\Psi$ of
$TS(\cA)$ with w*-topology into $[0,1]\times[0,1]\times\cM([0,1])$ with product
topology by $\Psi(\tau):=(\tau(e),\tau(f),\tilde\mu)$ where $\tilde\mu$ is the
distribution measure of $efe$ with respect to $\tau$. For each
$\tau=(\{\alpha_{ij}\}_{i,j=0}^1,\mu)\in TS(\cA)$ let $p:=\pi_\tau(e)$ and
$q:=\pi_\tau(f)$ in $(\pi_\tau(\cA)'',\tilde\tau)$, and let $e_{pqp}(\cdot)$ be
the spectral measure of $pqp$. From the structure theorem for two projections,
we get
\begin{eqnarray*}
\tilde\mu(\{0\})&=&\tilde\tau(e_{pqp}(\{0\}) \\
&=&{1\over2}\tilde\tau(\b1-
p\wedge q-p\wedge q^\perp-p^\perp\wedge q-p^\perp\wedge q^\perp) \\
&&\qquad+\tilde\tau(p\wedge q^\perp+p^\perp\wedge q+p^\perp\wedge q^\perp) \\
&=&{1\over2}(1-\alpha_{11}+\alpha_{10}+\alpha_{01}+\alpha_{00})
\end{eqnarray*}
and
$$
\tilde\mu(\{1\})=\tilde\tau(e_{pqp}(\{1\})=\tilde\tau(p\wedge q)=\alpha_{11}
$$
thanks to \eqref{F-3.1}. Hence it is straightforward to check that $\tau$
satisfies \eqref{F-3.2} if and only the following hold:
$$
\cases
\tau(e)=\alpha, \\
\tau(f)=\beta, \\
\tilde\mu(\{0\})=1-\min\{\alpha,\beta\}, \\
\tilde\mu(\{1\})=\max\{\alpha+\beta-1,0\}.
\endcases
$$
Furthermore, in this case we obviously have
$$
\tilde\mu=(1-\min\{\alpha,\beta\})\delta_0+\max\{\alpha+\beta-1,0\}\delta_1
+\rho\mu,
$$
where
\begin{equation}\label{F-3.4}
\rho=\min\{\alpha,\beta,1-\alpha,1-\beta\}
={1\over2}\Biggl(1-\sum_{i.j=0}^1\alpha_{ij}\Biggr).
\end{equation}
Based on Theorem \ref{T-2.2} together with these facts, to show the theorem,
it suffices to prove the following assertions:
\begin{itemize}
\item[(i)] If $\tau\in TS(\cA)$ and $(\tau(e),\tau(f))\ne(\alpha,\beta)$, then
$$
\inf_G\biggl[\limsup_{N\to\infty}{1\over N^2}\log\nu_N(G)\biggr]=-\infty.
$$
\item[(ii)] If $\tau\in TS(\cA)$ and $\Psi(\tau)=(\alpha,\beta,\tilde\mu)$, then
$$
\inf_G\biggl[\limsup_{N\to\infty}{1\over N^2}\log\nu_N(G)\biggr]
\le-\tilde I(\tilde\mu),
$$
$$
\inf_G\biggl[\liminf_{N\to\infty}{1\over N^2}\log\nu_N(G)\biggr]
\ge-\tilde I(\tilde\mu),
$$
where $\tilde I(\tilde\mu)$ is the rate function in Theorem \ref{T-2.2} and $G$
runs over neighborhoods of $\tau$.
\end{itemize}

When $(\tau(e),\tau(f))\ne(\alpha,\beta)$, choose $\eps>0$ such that
$\eps<|\tau(e)-\alpha|$ (or $\eps<|\tau(f)-\beta|$), and set
$G:=\{\tau'\in TS(\cA):|\tau'(e)-\alpha|<\eps
\ ({\rm or}\ |\tau'(f)-\beta|<\eps)\}$. Since
$\tau_N(e)=\tr_N(P(N))=k(N)/N\to\alpha$ and
$\tau_N(f)=\tr_N(Q(N))=l(N)/N\to\beta$ as $N\to\infty$, we get $\nu_N(G)=0$ for
large $N$ so that (i) follows.

To prove (ii), assume that $\Psi(\tau)=(\alpha,\beta,\tilde\mu)$. For any
neighborhood $\widetilde G$ of $\tilde\mu$, note that
$\Psi^{-1}([0,1]\times[0,1]\times\widetilde G)$ is a neighborhood of $\tau$ and
\begin{eqnarray*}
\nu_N\bigl(\Psi^{-1}([0,1]\times[0,1]\times\widetilde G)\bigr)
&=&{\rm Prob}\bigl(\{\Psi(\tau_N)
\in[0,1]\times[0,1]\times\widetilde G\}\bigr) \\
&=&{\rm Prob}\bigl(\{\widetilde R_N\in\widetilde G\}\bigr)
=\widetilde P_N(\widetilde G),
\end{eqnarray*}
where $\widetilde R_N$ is the empirical eigenvalue distribution of $P(N)Q(N)P(N)$
and $\widetilde P_N$ is its distribution on $\cM([0,1])$ (see the proof of
Theorem \ref{T-2.2}). Hence we have
$$
\inf_G\biggl[\limsup_{N\to\infty}{1\over N^2}\log\nu_N(G)\biggr]
\le\inf_{\widetilde G}\biggl[\limsup_{N\to\infty}{1\over N^2}
\log\widetilde P_N(\widetilde G)\biggr]\le-\tilde I(\tilde\mu)
$$
by Theorem \ref{T-2.2}. On the other hand, for any neighborhood $G$ of $\tau$,
one can choose $\eps>0$ and a neighborhood $\widetilde G$ of $\tilde\mu$ such
that $\Psi^{-1}\bigl((\alpha-\eps,\alpha+\eps)\times(\beta-\eps,\beta+\eps)
\times\widetilde G\bigr)\subset G$, which implies that
\begin{eqnarray*}
&&\liminf_{N\to\infty}{1\over N^2}\log\nu_N(G) \\
&&\quad\ge\liminf_{N\to\infty}{1\over N^2}\log
\nu_N\bigl(\Psi^{-1}\bigl((\alpha-\eps,\alpha+\eps)
\times(\beta-\eps,\beta+\eps)\times\widetilde G\bigr)\bigr) \\
&&\quad=\liminf_{N\to\infty}{1\over N^2}\log
{\rm Prob}\bigl(\bigl\{|\tr_N(P(N))-\alpha|<\eps,\,|\tr_N(Q(N))-\beta|<\eps,
\,\widetilde R_N\in\widetilde G\bigr\}\bigr).
\end{eqnarray*}
Since $|\tr_N(P(N))-\alpha|<\eps$ and $|\tr_N(Q(N))-\beta|<\eps$ for large $N$
(as in the proof of (i)), we have
$$
\liminf_{N\to\infty}{1\over N^2}\log\nu_N(G)
\ge\liminf_{N\to\infty}{1\over N^2}\log\widetilde P_N(\widetilde G)
\ge-\tilde I(\tilde\mu)
$$
by Theorem \ref{T-2.2}, and hence (ii) is proven. Finally, Theorem \ref{T-2.2}
proves the assertion on the minimizer as well (or this is a direct consequence
of the asymptotic freeness of $(P(N),Q(N))$).\qed

For $N\in\bN$ and $k\in\{0,1,\dots,N\}$ let $\cP(N,k)$ denote the set of all
$N\times N$ orthogonal projection matrices of rank $k$, and $\gamma_{N,k}$ be
the unitarily invariant measure on $\cP(N,k)$. We note that $\cP(N,k)$ is
identified with the homogeneous space $U(N)/(U(k)\oplus U(N-k))$ (or the
{\it Grassmannian manifold} $G(N,k)$) and $\gamma_{N,k}$ corresponds to the
measure on that space induced from the Haar probability measure on the unitary
group $U(N)$. In fact, an $N\times N$ unitarily invariant random projection
matrix of rank $k$ we have treated is standardly realized by $P\in\cP(N,k)$
distributed under $\gamma_{N,k}$.

Let $(p, q)$ be a pair of projections in a tracial
$W^*$-probability space $(\cM,\tau)$ and let $\alpha:=\tau(p)$ and
$\beta:=\tau(q)$. The {\it free entropy} $\chi(p,q)$ of $(p,q)$
proposed in \cite[14.2]{V6} by Voiculescu is defined as follows: Choose
sequences $k(N)$ and $l(N)$ such that $k(N)/N\to\alpha$ and $l(N)/N\to \beta$ 
as $N\to\infty$. For each $m\in\bN$ and $\eps>0$ set
\begin{eqnarray*}
&&\Gamma\bigl(p, q;k(N),l(N);N,m,\eps\bigr) \\
&&\qquad:=\Bigl\{(P,Q)\in \cP(N,k(N))\times \cP(N,l(N)):
\big|\tr_N(P_1\cdots P_m)-\tau(p_1\cdots p_m)\big|<\eps \\
&&\hskip6cm\mbox{for all $(P_j,p_j)\in \{(P,p),(Q,q)\}$, $1\le j\le m$}\Bigr\},
\end{eqnarray*}
and define
\begin{equation}\label{F-3.5}
\chi(p,q) 
:=\lim_{m\to\infty\atop\eps\searrow0}\limsup_{N\to\infty} 
{1\over N^2}\log
\bigl(\gamma_{N,k(N)}\otimes \gamma_{N,l(N)}\bigr)
\Bigl(\Gamma\bigl(p,q;k(N),l(N);N,m,\eps\bigr)\Bigr).
\end{equation}

Let $\cA$ be the $C^*$-algebra with two projection generators $e,f$ introduced
in the previous section. The {\it free entropy} of a tracial state
$\tau\in TS(\cA)$ is defined as $\chi(\pi_\tau(e),\pi_\tau(f))$ in the
tracial $W^*$-probability space $(\pi_\tau(\cA)'',\tilde\tau)$ obtained via the
GNS construction associated with $\tau$.

Next we identify the rate function in Theorem \ref{T-3.2} as the free entropy
$\chi(\tau)$ (up to a sign).

\begin{prop}\label{P-3.3}
The rate function in Theorem \ref{T-3.2} given for $\alpha=\tau(e)$ and
$\beta=\tau(f)$ is
$$
\cI(\tau)=-\chi(\tau).
$$
Moreover $\limsup$ can be replaced by $\lim$ in definition \eqref{F-3.5}.
\end{prop}

\proof
Let $p:=\pi_\tau(e)$, $q:=\pi_\tau(f)$ and $\tilde\mu$ be the distribution of
$efe$ with respect to $\tau$. In view of the form \eqref{F-3.3} of joint moments
of $e,f$ and the choices of $k(N),l(N)$ as above, one can easily see that for
each $m\in\bN$ and $\eps>0$
\begin{eqnarray*}
&&\Gamma(p,q;k(N),l(N);N,2m,\eps) \\
&&\quad=\Bigl\{(P,Q)\in\cP(N,k(N))\times\cP(N,l(N)):
\big|\tr_N((PQP)^k)-\tau((efe)^k)\big|<\eps,\,1\le k\le m\Bigr\}
\end{eqnarray*}
whenever $N$ is large enough. This implies that
$$
\bigl(\gamma_{N,k(N)} \otimes\gamma_{N,l(N)}\bigr)
\Bigl(\Gamma\bigl(p,q;k(N),l(N);N,2m,\eps\bigr)\Bigr)
=\widetilde P_N(\widetilde G(m,\eps)),
$$
where $\widetilde P_N$ is the distribution on $\cM([0,1])$ mentioned in the
proof of Theorem \ref{T-3.2} and $\widetilde G(m,\eps)$ is a neighborhood of
$\tilde\mu$ given by
$$
\widetilde G(m,\eps):=\biggl\{\tilde\mu'\in\cM([0,1]):
\bigg|\int x^k\,d\tilde\mu'(x)-\int x^k\,d\tilde\mu(x)\bigg|<\eps,\,
1\le k\le m\biggr\}.
$$
Now, as in the proof of \cite[5.6.2]{HP} we have the limit
\begin{eqnarray*}
&&\lim_{N\to\infty}{1\over N^2}\log
\bigl(\gamma_{N,k(N)}\otimes\gamma_{N,l(N)}\bigr)
\Bigl(\Gamma\bigl(p,q;k(N),l(N);N,2m,\eps\bigr)\Bigr) \\
&&\qquad=\lim_{N\to\infty}{1\over N^2}\log
\widetilde P_N(\widetilde G(m,\eps)),
\end{eqnarray*}
and the conclusion follows from Theorem \ref{T-3.2} and its proof.\qed

\bigskip
Theorem \ref{T-3.2} implies that the free entropy $\chi(p,q)$ of two projections
$p,q$ admits a maximal value, i.e., $\chi(p,q)=0$ if and only if $p,q$ are free.
Moreover, note by Proposition \ref{P-3.3} that the definition \eqref{F-3.5} of
$\chi(p,q)$ is independent of the choices of sequences $k(N)$ and $l(N)$, but
this fact is easy to directly verify. 

A further study of the free entropy $\chi(p_1,\dots,p_n)$ for general $n$-tuples
of projections as well as some related topics will be in a forthcoming paper
\cite{HU}.

\bigskip
\section{Applications of the contraction principle}
\setcounter{equation}{0}

Let $(P(N),Q(N))$ be as before, and let $\cA$ be the $C^*$-algebra of two
projection generators introduced in the previous section. Our large deviation in
Theorem \ref{T-3.2} is formulated on the tracial state space of $\cA$. The aim
of this section is to exemplify how Theorem \ref{T-3.2} implies, via the
contraction principle, the large deviation for the empirical eigenvalue density
of various random matrices made from $(P(N),Q(N))$.

For each selfadjoint element $h\in\cA$ and $\tau\in TS(\cA)$, let
$\lambda_h(\tau)$ denote the distribution measure of $h$ with respect to
$\tau$. Fixing $h$ we then have a map $\lambda_h:TS(\cA)\to\cM(\bR)$; in fact,
$\lambda_h(\tau)\in\cM([-\|h\|,\|h\|])$ for every $\tau\in TS(\cA)$. It is
straightforward to see that $\lambda_h$ is continuous with respect to
w*-topology on $TS(\cA)$ and weak topology on $\cM(\bR)$. Let
$\tau_N:=\tau_{P(N),Q(N)}$ be the random tracial state on $\cA$ induced by
$(P(N),Q(N))$ and $\nu_N$ the distribution on $TS(\cA)$ of $\tau_N$ (see
Section 3). We then notice that
$$
\nu_N\circ\lambda_h^{-1}(\Lambda)
={\rm Prob}(\{\tau_N\in\lambda_h^{-1}(\Lambda)\})
={\rm Prob}(\{\lambda_h(\tau_N)\in\Lambda\})
$$
for Borel sets $\Lambda\subset\cM(\bR)$. Since
$$
\int x^m\,d\lambda_h(\tau_N)(x)=\tau_N(h^m)
=\tr_N(h(P(N),Q(N))^m),\qquad m\in\bN,
$$
it follows that $\lambda_h(\tau_N)$ is nothing but the empirical eigenvalue
distribution of an $N\times N$ selfadjoint random matrix $h(P(N),Q(N))$ (via
``noncommutative functional calculus" mentioned in Section 3). Therefore, by
the {\it contraction principle} (see \cite[4.2.1]{DZ}), Theorem \ref{T-3.2}
implies the following:

\begin{thm}\label{T-4.1}
For every selfadjoint element $h\in\cA$, the empirical eigenvalue distribution
of $h(P(N),Q(N))$ satisfies the large deviation principle in the scale $1/N^2$
with the good rate function
$$
I_h(\mu):=\inf\{\cI (\tau):
\tau\in TS(\cA),\,\lambda_h(\tau)=\mu\}
$$
for $\mu\in\cM(\bR)$, and $\mu_0:=\lambda_h(\tau_0)$ is a unique minimizer of
$I_h$, where $\cI $ and $\tau_0$ are as in
Theorem \ref{T-3.2}.
\end{thm}

\begin{remark}\label{R-4.2}{\rm
For any unitary $u\in\cA$ define a map $\lambda_u:TS(\cA)\to\cM(\bT)$, $\bT$
being the unit circle, by letting $\lambda_u(\tau)$ the distribution of $u$
with respect to $\tau$. Then a similar large deviation is satisfied for the
empirical eigenvalue distribution of the unitary random matrix $u(P(N),Q(N))$
and the rate function $I_u$ is given in the same way as in Theorem \ref{T-4.1}.
}\end{remark}

In this way, for concrete applications, it remains only to find an explicit form
of the rate function $I_h$ (or $I_u$) as well as
that of the minimizer $\mu_0$. We present a few examples in the rest of the
section.

\bigskip\noindent
{\bf Example 4.3.}\enspace
Consider $h=ef+fe\in\cA$ and let $\tau=(\{\alpha_{ij}\}_{i,j=0}^1,\mu)\in
TS(\cA)$ as in Section 3. Since $e(t)f(t)+f(t)e(t)$ has the eigenvalues
$t\pm\sqrt t$, we get
\begin{eqnarray*}
\tau(\ffi(ef+fe))&=&(\alpha_{10}+\alpha_{01}+\alpha_{00})\ffi(0)
+\alpha_{11}\ffi(2) \\
&&\quad+\Biggl(1-\sum_{i,j=0}^1\alpha_{ij}\Biggr)
\int_0^1{\ffi(t+\sqrt t)+\ffi(t-\sqrt t)\over2}\,d\mu(t)
\end{eqnarray*}
for every continuous function $\ffi$ on $\bR$. By this expression and
\eqref{F-3.4}, whenever $\tau$ satisfies \eqref{F-3.2}, we have
\begin{eqnarray}\label{F-4.1}
\lambda_{ef+fe}(\tau)&=&\max\{|\alpha-\beta|,1-2\alpha,1-2\beta\}\delta_0
+\max\{\alpha+\beta-1\}\delta_2 \nonumber\\
&&\qquad+\rho(\mu\circ S^{-1}+\mu\circ T^{-1}),
\end{eqnarray}
where $S:(0,1)\to(0,2)$ and $T:(0,1)\to[-1/4,0)$ are given by
$St:=t+t\sqrt t$ and $Tt:=t-\sqrt t$. Hence the empirical eigenvalue
distribution of $P(N)Q(N)+Q(N)P(N)$ satisfies the large deviation in the scale
$1/N^2$ and the good rate function $\tilde I(\tilde\mu)$ for
$\tilde\mu\in\cM(\bR)$ is given by \eqref{F-2.9} if $\tilde\mu$ is of the form
in the right-hand side of \eqref{F-4.1} with $\mu\in\cM((0,1))$; otherwise
$\tilde I(\tilde\mu)=+\infty$. The minimizer of $\tilde I(\tilde\mu)$ is the
right-hand side of \eqref{F-4.1} with $\mu=\mu_0$, where $\rho\mu_0$ is the
continuous part of the measure \eqref{F-2.10}.

\bigskip\noindent
{\bf Example 4.4.}\enspace
Consider $h=ae+bf$ with $a,b\in\bR\setminus\{0\}$. Since $ae(t)+bf(t)$ has the
eigenvalues ${1\over2}(a+b\pm\sqrt{(a-b)^2+4abt})$, we get
\begin{eqnarray*}
\tau(\ffi(ae+bf))&=&\alpha_{00}\ffi(0)+\alpha_{10}\ffi(a)
+\alpha_{01}\ffi(b)+\alpha_{11}\ffi(a+b) \\
&&+\Biggl(1-\sum_{i,j=0}^1\alpha_{ij}\Biggr)\int_0^1{1\over2}
\biggl(\ffi\biggl({a+b-\sqrt{(a-b)^2+4abt}\over2}\biggr) \\
&&\hskip4cm+\ffi\biggl({a+b+\sqrt{(a-b)^2-4abt}\over2}\biggr)\biggr)\,d\mu(t)
\end{eqnarray*}
for every continuous function $\ffi$ on $\bR$ and
$\tau=(\{\alpha_{ij}\}_{i,j=0}^1,\mu)\in TS(\cA)$. Let $A,B$ be the first two
smallest numbers of $0,a,b,a+b$, and define $S:(0,1)\to(A,B)$ and
$T:(0,1)\to(a+b-B,a+b-A)$ by
$$
St:={a+b-\sqrt{(a-b)^2+4abt}\over2},\quad
Tt:={a+b+\sqrt{(a-b)^2+4abt}\over2}.
$$
When $\tau$ satisfies \eqref{F-3.2}, the above expression shows that
\begin{eqnarray*}
\lambda_{ae+bf}(\tau)&=&\max\{1-\alpha-\beta,0\}\delta_0
+\max\{\alpha-\beta,0\}\delta_a \nonumber\\
&&\quad+\max\{\beta-\alpha,0\}\delta_b
+\max\{\alpha+\beta-1,0\}\delta_{a+b} \nonumber\\
&&\qquad+\rho(\mu\circ S^{-1}+\mu\circ T^{-1}).
\end{eqnarray*}
Hence the empirical eigenvalue distribution of $aP(N)+bQ(N)$ satisfies the
large deviation and the good rate function as well as its minimizer is
determined similarly to the above example.

Let us express the rate function $\tilde I(\tilde\mu)$ and the minimizer
$\tilde\mu_0$ more explicitly. When $\mu\in\cM((0,1))$, the measure
$\nu:={1\over2}(\mu\circ S^{-1}+\mu\circ T^{-1})$ is supported in
$(A,B)\cup(a+b-B,a+b-A)$ and symmetric at $(a+b)/2$ so that
$\mu=2\nu\circ S|_{(A,B)}=2\nu\circ T|_{(a+b-B,a+b-A)}$. Since $St=x$ (or
$Tt=x$) implies $t=(x-a)(x-b)/ab$, we get
$$
\int_0^1\log t\,d\mu(t)=2\int_A^B\log{(x-a)(x-b)\over ab}\,d\nu(x)
=2\int_{a+b-B}^{a+b-A}\log{(x-a)(x-b)\over ab}\,d\nu(x)
$$
so that
$$\int_0^1\log t\,d\mu(t)
=\int_{(A,B)\cup(a+b-B,a+b-A)}\log{(x-a)(x-b)\over ab}\,d\nu(x).
$$
Similarly,
$$
\int_0^1\log(1-t)\,d\mu(t)
=\int_{(A,B)\cup(a+b-B,a+b-A)}\log{x(a+b-x)\over ab}\,d\nu(x).
$$
On the other hand, we get
\begin{eqnarray*}
\Sigma(\mu)&=&4\int_A^B\int_A^B\log
\bigg|{(x-a)(x-b)\over ab}-{(y-a)(y-b)\over ab}\bigg|\,d\nu(x)\,d\nu(y) \\
&=&4\int_A^B\int_A^B\log
\bigg|{(x-a)(a+b-x-y)\over ab}\bigg|\,d\nu(x)\,d\nu(y) \\
&=&2\Sigma(\nu)-\log|ab|.
\end{eqnarray*}
Consequently, the rate function $\tilde I(\tilde\mu)$ is written as
\begin{eqnarray*}
\tilde I(\tilde\mu)&=&-2\rho^2\Sigma(\nu)
-\rho|\alpha-\beta|\int_{(A,B)\cup(a+b-B,a+b-A)}\log|(x-a)(x-b)|\,d\nu(x) \\
&&\qquad-\rho|\alpha+\beta-1|
\int_{(A,B)\cup(a+b-B,a+b-A)}\log|x(a+b-x)|\,d\nu(x) \\
&&\qquad+C +\rho\max\{\alpha,\beta,1-\alpha,1-\beta\}\log|ab|
\end{eqnarray*}
if $\tilde\mu\in\cM(\bR)$ is of the form
\begin{eqnarray*}
\tilde\mu&=&\max\{1-\alpha-\beta,0\}\delta_0
+\max\{\alpha-\beta,0\}\delta_a \\
&&\quad+\max\{\beta-\alpha,0\}\delta_b
+\max\{\alpha+\beta-1,0\}\delta_{a+b}+2\rho\nu
\end{eqnarray*}
with $\nu\in\cM((A,B)\cup(a+b-B,a+b-A))$ symmetric at $(a+b)/2$; otherwise
$\tilde I(\tilde\mu)=+\infty$.

Moreover, by transforming the continuous part of \eqref{F-2.10}, the explicit
form of the minimizer $\tilde\mu_0$ can be easily computed as follows:
\begin{eqnarray}\label{F-4.2}
\tilde\mu_0&=&\max\{1-\alpha-\beta,0\}\delta_0
+\max\{\alpha-\beta,0\}\delta_a \nonumber\\
&&\quad+\max\{\beta-\alpha,0\}\delta_b
+\max\{\alpha+\beta-1,0\}\delta_{a+b} \nonumber\\
&&+{\big|x-{a+b\over2}\big|\sqrt{-(x-A_0)(x-B_0)(x-a-b+B_0)(x-a-b+A_0)}
\over\pi|x(x-a)(x-b)(x-a-b)|} \nonumber\\
&&\hskip5cm\times\b1_{(A_0,B_0)\cup(a+b-B_0,a+b-A_0)}(x)\,dx,
\end{eqnarray}
where
$$
A_0:={a+b-\sqrt{(a-b)^2+4ab\eta}\over2},\quad
B_0:={a+b-\sqrt{(a-b)^2+4ab\xi}\over2}
$$
(or exchange $A_0,B_0$ depending on the sign of $ab$) with $\xi,\eta$ in
\eqref{F-2.11}. As is guaranteed by the asymptotic freeness (\cite{Vo}) of
$(P(N),Q(N))$, the minimizer $\tilde\mu_0$ is equal to the distribution of
$ap+bq$ where $(p,q)$ is a pair of free projections in a tracial
$W^*$-probability space $(\cM,\tau)$ with $\tau(p)=\alpha$ and $\tau(q)=\beta$.
In fact, the distribution was computed in \cite{AY} by use of $R$-transform.

Although one can prove the large deviation result for the empirical eigenvalue
density of $aP(N)+bQ(N)$ (also $P(N)Q(N)+Q(N)P(N)$) based on the joint
eigenvalue distributions given in Corollary \ref{C-1.3}, our stress is that this
is just a particular case of grand Theorem \ref{T-4.1} (or Theorem \ref{T-3.2}). 

\bigskip\noindent
{\bf Example 4.5.}\enspace
For unitaries we consider a simple example $u=e^{\im\pi e}e^{-\im\pi f}$. Since
the eigenvalues of $e^{\im\pi e(t)}e^{-\im\pi f(t)}$ are
$2t-1\pm2\im\sqrt{t(1-t)}=e^{\pm\im\theta(t)}$ where
$\theta(t):=\cos^{-1}(2t-1)$ for $t\in(0,1)$, we get
\begin{eqnarray*}
\tau(\ffi(u))&=&(\alpha_{11}+\alpha_{00})\ffi(1)
+(\alpha_{10}+\alpha_{01})\ffi(-1) \\
&&\qquad+\Biggl(1-\sum_{i,j=0}^1\alpha_{ij}\Biggr)
\int_0^1{\ffi(e^{\im\theta(t)})+\ffi(e^{-\im\theta(t)})\over2}\,d\mu(t)
\end{eqnarray*}
for every continuous function $\ffi$ on $\bT$ and
$\tau=(\{\alpha_{ij}\}_{i,j=0}^1,\mu)\in TS(\cA)$. When $\tau$ satisfies
\eqref{F-3.2}, this implies that
$$
\lambda_u(\tau)=|\alpha+\beta-1|\delta_1+|\alpha-\beta|\delta_{-1}
+\rho(\mu\circ\theta^{-1}+\mu\circ\tilde\theta^{-1}),
$$
where $\tilde\theta(t):=-\theta(t)$ for $t\in(0,1)$. For $\mu\in\cM((0,1))$ let
$\nu:={1\over2}(\mu\circ\theta^{-1}+\mu\circ\tilde\theta^{-1})$, which is a
probability measure on $\bT$ symmetric for the real axis. We then have
$$
\int_0^1\log t\,d\mu(t)
=\int_\bT\log{1+\cos\theta\over2}\,d\nu(e^{\im\theta}),
$$
$$
\int_0^1\log(1-t)\,d\mu(t)
=\int_\bT\log{1-\cos\theta\over2}\,d\nu(e^{\im\theta}),
$$
$$
\Sigma(\mu)=\iint_{\!\!\!\bT^2}\log|\cos\theta-\cos\psi|
\,d\nu(e^{\im\theta})\,d\nu(e^{\im\psi})-\log2.
$$
Hence we see by Remark \ref{R-4.2} that the empirical eigenvalue distribution of
$e^{\im\pi P(N)}e^{-\im\pi Q(N)}$ satisfies the large deviation in the scale
$1/N^2$ and the rate function is given by
\begin{eqnarray*}
\tilde I(\tilde\mu)
&=&-\rho^2\iint_{\!\!\!\bT^2}\log|\cos\theta-\cos\psi|
\,d\nu(e^{\im\theta})\,d\nu(e^{\im\psi}) \\
&&-\rho|\alpha-\beta|\int_\bT\log(1+\cos\theta)\,d\nu(e^{\im\theta})
-\rho|\alpha+\beta-1|\int_\bT\log(1-\cos\theta)\,d\nu(e^{\im\theta}) \\
&&\qquad+C +\rho\max\{\alpha,\beta,1-\alpha,1-\beta\}\log2
\end{eqnarray*}
if $\tilde\mu\in\cM(\bT)$ is of the form
$\tilde\mu=|\alpha+\beta-1|\delta_1+|\alpha-\beta|\delta_{-1}+2\rho\nu$
with $\nu\in\cM(\bT)$ having no atoms at $\pm1$ and symmetric for the real axis;
otherwise $\tilde I(\tilde\nu)=+\infty$. The minimizer $\tilde\mu_0$ is also easy
to compute as
\begin{eqnarray}\label{F-4.3}
\tilde\mu_0&=&|\alpha+\beta-1|\delta_1+|\alpha-\beta|\delta_{-1} \nonumber\\
&&\quad+{\sqrt{-(\cos\theta+1-2\xi)(\cos\theta+1-2\eta)}\over|\sin\theta|}
\b1_{(\theta_1,\theta_2)\cup(-\theta_2,-\theta_1)}(\theta)
\,{d\theta\over2\pi},
\end{eqnarray}
where $\theta_1:=\cos^{-1}(2\eta-1)$ and $\theta_2:=\cos^{-1}(2\xi-1)$. This
measure is the distribution of $e^{\im\pi p}e^{-\im\pi q}$ for free projections
$p,q$ sometimes mentioned above. It may be natural that this distribution is
rather different (except the same atomic parts) from that of  $e^{\im\pi(p-q)}$
computed from \eqref{F-4.2}. In particular, when $\alpha=\beta=1/2$ so that
$\xi=0$ and $\eta=1$, the minimizer \eqref{F-4.3} is the uniform measure on
$\bT$ but \eqref{F-4.2} induces the arcsine law on the angular variable
$(-\pi,\pi)$.

\end{document}